\documentclass[a4paper,12pt]{amsart}
\usepackage{amsmath,amssymb,amsfonts,amsthm}

\usepackage[bottom]{footmisc}
\usepackage[utf8x]{inputenc}
\usepackage{tikz,relsize} 
\usepackage[cmtip,arrow]{xy}
\usepackage{pb-diagram,pb-xy}
\usepackage{fontenc}
\usepackage{hyperref}

\newtheorem{prop}{Proposition}
\newtheorem{rem}{Remark}
\newtheorem{exe}{Example}
\newtheorem{theo}{Theorem}
\newtheorem{Df}{Definition}
\newtheorem{lm}{Lemma}
\newtheorem{cor}{Corollary}
\usepackage[mathscr]{eucal}
\usepackage{indentfirst}
\usepackage{calligra}
\usepackage{eucal}
\usepackage[mathscr]{eucal}
\usepackage[english]{babel}
\usepackage{underscore}
\newcommand{\R} {{\mathbb R}}
\newcommand{\Q} {{\mathbb Q}}
\newcommand{\Z} {{\mathbb Z}}

\newcommand{\A} {{\mathbb A}}

\newcommand{\C} {{\mathbb C}}
\newcommand{\F} {{\mathbb F}}

\newcommand\blfootnote[1]{%
\begingroup
\renewcommand\thefootnote{}\footnote{#1}%
\addtocounter{footnote}{-1}%
\endgroup}
\begin{document}
\def\ns#1{\mathbb{#1}}
\title{Properties of the combinatorial Hantzsche-Wendt groups} 
\author{J. Popko and A. Szczepa\'nski}

\address{Jerzy Popko, Andrzej Szczepa\'nski\\
Institute of Mathematics,\\
University of Gda\'nsk,\\
ul. Wita Stwosza 57,\\
80-308 Gda\'nsk, Poland}
\email{jpopko@mat.ug.edu.pl, matas@ug.edu.pl}

\begin{abstract}The combinatorial Hantzsche-Wendt group $G_n =\\ \langle x_1,...,x_n\mid x_{i}^{-1}x_{j}^{2}x_{i}x_{j}^{2},\forall i\neq j\rangle$
	was defined by W. Craig and P. A. Linnell in \cite{CL}. For $n=2$ it is a fundamental group of 3-dimensional oriented
	flat manifold with non cyclic holonomy group. We calculate the Hilbert-Poincar\'e 
	series of $G_n, n\geq 1$ with $\Q$ and $\F_2$ coefficients. Moreover, we prove that the cohomological dimension of $G_n$
	is equal to $n+1.$ Some other properties of this group are also considered.
\end{abstract}

\blfootnote{2020 {\it Mathematics Subject Classification.} Primary: 20J05; Secondary: 20H15 20E06 55R20\\
{\it Key words and phrases.} Hantzsche-Wendt groups, Hilbert-Poincar\'e series, cohomological dimension\\
{\it Corresponding author:} A. Szczepa\'nski}



\date{\today}

\maketitle
\vskip 5mm
\noindent
\section{Introduction}
Let $\Gamma_3$ be the fundamental group of the oriented flat 3-manifold with non-cyclic holonomy,
which was the first time defined by W. Hantzsche $\&$ H. Wendt and W. Nowacki in 1934, see
\cite{HW}, \cite{N}. From \cite[ch. 9]{S}, $\Gamma_3$ is a torsion free crystallographic group of a rank 3.
Where, by crystallographic group of dimension $n$ we understand a discrete and cocompact subgroup
of the group $E(n) = O(n)\ltimes\R^n$ of isometries of the Euclidean space $\R^n.$ From the Bieberbach theorems
\cite{S} any crystallographic group $\Gamma$ of rank $n$ defines a short exact sequence
\begin{equation}\label{short} 
1\to\Z^n\to\Gamma\to H\to 1,
\end{equation}
where $\Z^n$ is the free abelian subgroup of all translations of $\Gamma$ and $H$ is a finite group, called the
holonomy group of $\Gamma.$ In the case of $\Gamma_3$ the group $H =\Z_2\oplus\Z_2.$ As a subgroup of $E(3)$
	{\scriptsize $$\Gamma_3 = \text{gen}\{A =(\left[
	\begin{array}{clcr}
	1 & 0 & 0\\
	0 & -1 & 0\\
	0 &  0 & -1
	\end{array}
	\right],(1/2,1/2,0)),
	B = (\left[
	\begin{array}{clcr}
	-1 & 0 & 0\\
	0 & 1 & 0\\
	0 & 0 & -1
	\end{array}
	\right],(0,1/2,1/2))\}.
	$$}
The Hantzsche-Wendt groups/manifolds are also defined in higher odd dimensions, as fundamental groups of oriented
flat manifolds of dimensions, $n\geq 3$ with holonomy group $(\Z_2)^{n-1}.$ We shall denote them by $\Gamma_n$, see \cite[ch. 9]{S}.
From the Bieberbach theorems there exist, for given $n,$ a finite number $L(n)$ of Hantzsche-Wendt groups (HW groups),
up to isomorphism. However, the number $L(n)$ growths exponentially, see \cite[Theorem 2.8]{MR}. Let us define an
example of the HW group $\Gamma_n$ of dimension $\geq 3$ which is a generalization of $\Gamma_3.$
\begin{exe} Let $n$ be an odd number. Then $\Gamma_n =$
{\scriptsize $$= \text{gen}\{\gamma_i=(\left[
\begin{array}{clcrclcrcl}
-1 & 0 & ...  & 0 & 0 & 0 & ... & 0 & 0 & \\
0  &-1 & ...  & 0 & 0 & 0 & ... & 0 & 0 & \\
... & \\
0 & 0 & ... & -1 & 0 & 0 & ... & 0 & 0 \\
0 & 0 & ... & 0 &  1 & 0 & ... & 0 & 0 \\
0 & 0 & ... & 0 & 0 & -1 & ... & 0 & 0 \\
... \\
0 & 0 & ... & 0 & 0 & 0 & ... & -1 & 0 \\
0 & 0 & ... & 0 & 0 & 0 & ... & 0 & -1 
\end{array}\right],(0,...,0,1/2,1/2,0,...,0))\},$$}
where $1$ is at the $i$-th place and the first $1/2$ is at the $i$-th place, $1\leq i\leq n-1.$
\end{exe}

In 1982, see \cite{S}, the second author proved that for odd $n\geq 3,$ 
the manifolds $\R^n/\Gamma_n$ are rational homology spheres. Moreover,
for $n\geq 5$ the
commutator subgroup of the group $\Gamma_n$ is equal to the translation
subgroup \cite[Theorem 9.3]{S}.
Moreover, for $m\geq 7$ there exist many isospectral HW-manifolds 
non pairwise homeomorphic, \cite[Corollary 3.6]{MR}. 
HW groups have an interesting connection with Fibonacci groups
(see below) and the theory of quadratic forms over the field $\F_2,$
\cite[Theorem 9.5]{S}. HW-manifolds have no $\operatorname{Spin}$ or $\operatorname{Spin^{\C}}$-structures, 
\cite{L} and \cite[page 109]{S}. 
Finally HW manifolds are cohomological rigid that means two HW manifolds are homeomorphic if and only if their
cohomology rings over $\F_2$ are isomorphic, \cite{PS}.
\vskip 1mm
$G$ is called a unique product group if given two nonempty finite subset $X,Y$ of $G,$
then exists at least one element $g\in G$ which has a unique representation $g=xy$ with $x\in X$ and $y\in Y.$
A unique product group is torsion free, though the converse is not true in general.
The original motivation for studying unique product groups was the Kaplansky zero divisor conjecture, namely that 
if $k$ is a field and $G$ is a torsion free group, then $kG$ is a domain.
It was proved in 1988 \cite{Prom} that the group $G_2$ is a nonunique product group.
To prove it the author uses the combinatorial presentation (\cite[Lemma 13.3.1, p. 606-607]{Pas})
\begin{equation}\label{pres}
\Gamma_3 = \langle x,y\mid x^{-1}y^{2}xy^{2}, y^{-1}x^{2}yx^{2}\rangle.
\end{equation}
However the counterexample to the Kaplansky unit conjecture was given in 2021 by G. Gardam \cite{Ga}.
Again the counterexample was found in the group ring $\F_2[\Gamma_3].$
The Kaplansky unit conjecture states that every unit in $K[G]$ is of the form $kg$ for
$k ∈ K \setminus \{0\}$ and $g \in G.$
\vskip 1mm
\noindent
In \cite{CL} the following generalization of $\Gamma_3$ is proposed.
\begin{Df} By a combinatorial Hantzsche-Wendt group we shall understand a finitely presented group
$$G_n = \langle x_1, . . . , x_n | x_{i}^{-1}x_{j}^{2}x_{i}x_{j}^{2}\hskip 3mm \forall\hskip 2mm i\neq j\rangle.$$
\end{Df}
\noindent
It is easy to see that, $G_{0} = 1$ and $G_1 =\Z.$ Moreover $G_2$ is the Hantzsche - Wendt group of dimension 3.

Let 
\begin{equation}\label{An}
\Z^n\simeq\A_n\simeq\langle x_1^2, x_2^2,...,x_n^2\rangle,
\end{equation}	
be a free abelian subgroup of $G_n.$
In \cite[Lemma 3.1]{CL} is proved that $\Z^{n}\lhd G_n.$ 
Later we shall denote $\A_n$ by $\Z^n.$
Moreover, $W_n = G_n/\Z^n =\langle x_1,...,x_n \mid  x_{1}^{2}, x_{2}^{2},...,x_{n}^{2}\rangle \simeq \ast^{n}\hskip 2mm\Z_2$. 
Finally in \cite[Theorem 3.3]{CL} it is proved that $G_n$ is torsion free for all $n\geq 1.$
This is also the corollary from Theorem \ref{theo2}.
For any  $1\leq m\leq n,$ $G_m$ embeds in $G_n$
and for $n\geq 2, G_n$ is a nonunique product group \cite[Corollary 3.5]{CL}.
Another interesting result of \cite[Theorem 3.6]{CL} is the following.
There is for $n\geq 3$ and odd a surjective homomorphism $\Phi_{n}:G_{n-1}\to\Gamma_n.$ It is easy to see
that $\Phi_{n}(\Z^{n-1})$ is a free abelian subgroup of the translation subgroup of $\Gamma_n$ of a rank $n-1.$
Since $\Gamma_n/\Phi_{n}(\Z^{n-1})$ is an infinite group and  $Ker(\Phi_{n})\cap \Z^{n-1} = {1}$
then $Ker(\Phi_{n})$ is an infinitely generated free group. (See \cite[Theorem 3.6]{CL} and \cite[page 87]{AT}.)

At that point we would like to mention the following related result, see \cite{LM}.
Recall that the Fibonacci group F(r, n) is defined by the presentation
$$F(r, n) = \langle a_0, . . . , a_{n−1} | a_{i}a_{i+1}\cdot \cdot \cdot a_{i+r−1} = a_{i+r} , 0\leq i\leq n − 1\rangle,$$
where the indices are understood modulo $n.$ There exists a connection of these groups with our family $G_n.$
We know that $F(2, 6)$ is isomorphic to $\Gamma_3$, and there is, for any $n\geq 3$ an epimorphism $\Psi_n:F(n-1,2n)\to\Gamma_n.$
\vskip 5mm
\noindent
In the first part of a paper we shall show two models of $BG_n$ (or $K(G_n,1)).$
They are a topological realization of two algebraic representations of $G_n.$ The first model is an appriopriate gluing
of $n$ copies of generalized fat Klein bottles. It corresponds to an isomorphism of $G_n$ with $\ast_{\Z^n}^{n}K_n$ where $K_n$ is
a generalized Klein bottle crystallographic group amalgamated over the translation lattices.
The second model is some Borel construction. It corresponds to the representation of $G_n$ as an extension:
\begin{equation}\label{maine}
1\to\Z^n\to G_{n}\to W_n\to 1.
\end{equation}
From the first model we obtain that the cohomological dimension of $G_n$ is equal to $n+1$ for $n > 1.$

In the second part we calculate the Hilbert-Poincar\'e series of $G_n, n\geq 1$ with $\Q$ and $\F_2$ coefficients
and explain the algebra structure of the cohomology. Here our main tools will be the Lyndon-Hochschild-Serre (LHS) spectral sequence
of the group extension (\ref{maine}). The case with $\F_2$ coefficients uses a multiplicative structure of LHS.
In the $\F_2$ case it is enough to use $E^{\star,\star}_{3}$ groups, but for rational coefficients we only need the 
$E^{\star,\star}_{2}-$ terms.
(See formulas \ref{Pp1} and \ref{Pp2}.)
\vskip 1mm
In the last part we calculate some other invariants and properties of $G_n.$ For example their abelianization and the Euler characteristic.
\vskip 5mm
\noindent
{\bf Acknowledgements} 
The authors are grateful to an anonymous referee for many helpful remarks and suggestions.
Moreover we would like to thank Giles Gardam for pointing out a small typo in an earlier version of the paper.
\section{Two models of $BG_n$}
\subsection{Gluing fat Klein bottles}
We start with an example.
\begin{exe} 
Let $K_{-}$ be the fundamental group of the Klein bottle and $\Z^2$ its maximal abelian subgroup of index two.	
It is well known (see \cite[Chapter 8.2, p. 153]{hill}) that $\Gamma_3\simeq K_{-}\ast_{\Z^2}K_{-}.$
\end{exe} 

A generalization of the above example gives us the followng characterization of the combinatorial Hantzsche - Wendt group.
Let $G_{n}^{(i)}$ denotes the subgroup of $G_n$ generated by $\{x_{i}\}$ and the abelian subgroup (\ref{An}) $\Z^n.$
We shall call it a generalized Klein bottle.
\begin{prop}	
The natural group homomorphism 
\begin{equation}\label{isom}
	\ast_{\Z^n}G_{n}^{(i)}\to G_n
\end{equation}
is an isomorphism.
\end{prop}
\vskip 1mm
\noindent
{\bf Proof:} 
This follows from the definition and the structure of the free product with amalgamation.
\vskip 3mm
\hskip 120 mm $\Box$
\vskip 2mm
\noindent
$G_{n}^{(i)}$ is a torsion free crystallographic group of dimension $n$ and acts freely on $\R^n$ (in a way analogous to $K_{-})$
so has a classifying space which is an $n$ dimensional closed flat manifold $K^{i}$ (the generalized Klein bottle).
A topological interpretation of the isomorphism (\ref{isom}) gives us a $n+1$ dimensional classifying space $BG_n$ as (homotopically)
gluing together $n$ generalized Klein bottles $K^{(1)}, K^{(2)},...,K^{(n)}$ along a common $n$ dimensional torus $\R^n/\Z^n.$
This space has dimension $n+1$ since we must convert maps $\R^n/\Z^n\to K^{(i)}$ to inclusions.
More precisely it may be done as follows. Let us define an action of $G_n$ on $\R^n$ by 
$$x_{i}(v)_{i} = v_i + 1/2\hskip 2mm\text and\hskip 2mm x_{i}(v)_{j} = -v_{j}, j\neq i,$$
where $v = (v_1,v_2,...,v_n)$ and an action on a segment $I = [-1,1]$ by $x_{i}(t) = -t, t\in I.$
\begin{Df}
By a fat Klein bottle we shall understand the space $B_{n}^{(i)}:= (\R^n\times I)/G_{n}^{(i)}.$
\end{Df}
Let $S_n: = \R^n/\Z^n$. Let us define maps $\alpha^{(i)}:S_n\to B_{n}^{(i)}$ by the formula $\alpha^{(i)}(v) = [(v,1)].$
\begin{Df}
By the space $B_n$ we shall understand a $colim$ of a diagram formed from maps $\alpha^{(i)},$
i.e. 
$$B_n:= colim_{i}\alpha^{(i)}.$$
\end{Df}
\begin{theo}
The above space $B_n$ is a classifying space for $G_n.$
\end{theo}
\vskip 1mm
\noindent
{\bf Proof:} From the definition the action of the subgroup $G_{n}^{(i)}$ on $\R^n$ is free and the orbit space
$K_{n}^{(i)}$ was called a generalized Klein bottle. Moreover, the fat Klein bottle is $(n+1)$
dimensional compact manifold with boundary and the projection on the first factor gives a bundle
$B_{n}^{(i)}\to K_{n}^{(i)}$ with fiber $I,$ hence in particular it is a homotopy equivalence.
Finally, the map $\alpha^{(i)}$ is an embedding on the boundary of $B_{n}^{(i)}.$
However, more geometricaly we may write $B_n:= \bigcup_{i }B_{n}^{(i)}$ treating the maps $\alpha^{(i)}$
as identifications (so $S_n\subset B_{n}^{(i)}$). In other words $B_n$ is obtained from $n$ copies
of a fat Klein bottle by an appriopriate identification of the boundaries of different copies.
To finish our proof we observe that $\pi_{1}(B_n)\simeq G_n$ after van Kampen theorem.
The space $B_n$ is aspherical after JHC Whitehead's theorem in \cite[Theorem 5]{Wh}.
\vskip 2mm
\hskip 120mm $\Box$
\begin{cor}
For $n > 1$ the cohomological dimension of $G_n$ is equal to $n+1.$
\end{cor}
\vskip 1mm
\noindent
{\bf Proof:} For brevity we write $B = B_n$ and $S = S_n.$ From the properties of $B$ we have, that cd$G_n\leq n+1.$
Let $H$ denote cohomology with $\F_2$ coefficients. We have an exact sequence
$$H^n(S)\to H^{n+1}(B,S)\to H^{n+1}(B).$$
Since dim$H^n(S) = 1$ and dim$H^{n+1}(B,S) = n,$ then dim$H^{n+1}(B)\geq n - 1.$ Hence cd$G_n\geq n+1$ for $n > 1.$
\vskip 2mm
\hskip 120mm $\Box$
\begin{rem}
1.- The space $B_n$ is for $n = 1$ a M\"obius band, for $n = 2$ a closed manifold (a classical
$3$-dimensional Hantzsche-Wendt manifold). However for $n > 2$ it is nonmanifold, since
there is singularity along $S_n.$
\end{rem}

\subsection{Borel construction (homotopy quotient)}
	
Let $G$ be a discrete group and let $p_G : EG\to BG$ be the universal $G$ bundle.
The assignment $G\mapsto p_G$ may be done functorial in the group $G$ and respecting products.
If $X$ is some $G$-space then the space
$$X_{G} := (EG\times X)/G$$
is called the Borel construction on $X,$ \cite[p. 10]{AP}. Here $G$ acts on $EG\times X$ diagonally.
Let $f_X :X_G\to BG$ be the quotient map. It is a fibration with fiber $X.$
It is easy to see that, if $X$ is aspherical then $X_G$ is also aspherical.

\begin{Df} (morphisms between maps)
If $f:X_1\to X_2$ and $g:Y_1\to Y_2$ then a morphism from $f$ to $g$
	is a pair $(m_1,m_2)$ where $m_i:X_i\to Y_i$ and $gm_1 = m_{2}f.$
\end{Df}

The operation of taking pullback along $\psi$ is denoted by $\psi^{\star}.$
We shall write $f\simeq m_{2}^{\star}(g)$ if $(m_1,m_2):f\to g$ and $m_1$ is an isomorphim on fibres.

Let $\xi,\eta\in W_2$ be generators of order 2 and $D :=\Z_2\oplus\Z_2.$
The abelianization of $W_2$ defines
\begin{equation}\label{dihed}
1\to\Z\simeq (\langle \xi\eta)^2\rangle \to W_{2}\stackrel{\alpha}\to D\to 1.
\end{equation}

Let $\Sigma$ be the unit circle on the complex plane.
Define an action of $D$ on $\Sigma$ by formulas
$$\xi(z) = \bar{z} \hskip 2mm \text{and}\hskip 2mm \eta(z) = -\bar{z}.$$
Denote a resulting $D$-space by $U.$ In the above language we have a map $f_{U}:U_{D}\to BD$
and and we can observe that $\pi_{1}(f_{U})\simeq\alpha.$

We define (for $i = 1,2,...,n$) homomorphisms $\phi_i:G_n\to W_2$
\begin{equation}\label{fii}
\phi_{i}(x_i) = \xi\eta\hskip 2mm\text{and}\hskip 2mm \phi_{i}(x_j) = \xi\hskip 2mm\text{for}\hskip 2mm j\neq i.
\end{equation}
	
Then $\phi:= (\phi_{1},...,\phi_n)$ gives a homomorphism from $G_n\to (W_2)^n.$
Let $q_{n} : G_n\to W_n$ be the canonical surjection. The homomorphism $\phi$ factorizes and we
obtain a map $(\phi,\psi):q_n\to\alpha^{n}$ and $\phi$ is an isomorphism on fibres. We have.
\begin{lm}
$q_n\simeq\psi^{\star}\alpha^{n}.$
\end{lm} 
\vskip 1mm
\noindent
{\bf Proof:}
With support of (\ref{dihed}) we have the following commutative diagram
$$
\begin{diagram}	
	\node{1}\arrow{e}\node{\Z^n}\arrow{s,r}{\simeq}\arrow{e,t}{i}\node{G_{n}}\arrow{s,r}{\phi}\arrow{e,t}{q_n}\node{W_{n}}\arrow{s,r}{\psi}\arrow{e}\node{1}\\
	\node{1}\arrow{e}\node{\Z^{n}}\arrow{e,t}{i_{1}}\node{(W_{2})^{n}}\arrow{e,t}{\alpha^{n}}\node{(\Z_{2}\oplus\Z_{2})^{n}}\arrow{e}\node{1}
\end{diagram},$$
\vskip 1mm
\centerline{Diagram 1}
\vskip 1mm
where $i, i_1$ are inclusions.  
\vskip 2mm
\hskip 120mm $\Box$

Define a $W_n$ action on the space $\Sigma^{n}$
\begin{equation}\label{action}
x_i(z)_i = - z_{i}\hskip 2mm\text{and}\hskip 2mm x_i(z)_j = \bar{z_j}\hskip 2mm\text{for}\hskip 2mm j\neq i,
\end{equation}
where $z = (z_1,z_2,...,z_n).$
\vskip 1mm
\noindent
Denote the resulting $W_n$-space by $T_n.$
\begin{prop}\label{prop1}
$$\pi_{1}(f_{T_n})\simeq q_n$$
in particular $\pi_{1}((T_n)_{W_n})\simeq G_n.$
\end{prop}
\vskip 2mm
\noindent
{\bf Proof:} The action on $T_n$ is obtained by composing the product $D^n$ action with the homomorphism $\psi$
(i.e. $w(x) = \psi(w)(z)$ for $w\in W_n$). Hence, from naturality we have a map of fibrations 
$$(\hat{\psi}_1,\hat{\psi}_2) = \hat{\psi}:f_{T_n}\to f_{U^n}.$$
Applying $\pi_1$ we get
$$\pi_1(\hat{\psi}):\pi_1(f_{T_n})\to \pi_1(f_{U^n}).$$
The map $\hat{\psi}$ gives an isomorphism (identity) on fibres so the map $\pi_1(\hat{\psi})$ also gives an isomorphism on fibres
by an application of the long exact sequence for fibrations. We have (for codomain components)
$(\hat{\psi})_2 = B\psi$ so $\pi(\hat{\psi})_2 = \pi(B\psi) = \psi.$ Hence
$$\pi_{1}(f_{T_n})\simeq\psi^{\star}\pi_1(f_{U^n}).$$
But $$\psi^{\star}\pi_1(f_{U^n})\simeq \psi^{\star}\pi_1(f_{U}))^n\simeq\psi^{\star}\alpha^n\simeq q_n.$$
\vskip 2mm
\hskip 120mm $\Box$
\begin{theo}\label{theo2}
$$(T_n)_{W_n} = K(G_n,1).$$
\end{theo}
\vskip 2mm
\noindent
{\bf Proof:} The space $(T_n)_{W_n}$ is aspherical because $T_n$ is. And it has the appriopriate fundamental group by proposition \ref{prop1}.
\vskip 2mm
\hskip120mm $\Box$

Let $B = \bigvee_{1}^{n}\R P(\infty).$ The space $B = K(W_n,1)$, cf. \cite{Wh}.
Let $E\to B$ be the universal- covering. Then
\begin{cor}
We have the fibration
\begin{equation}\label{fib}
T^n\to (T^n)_{W_n}\to E/W_n,
\end{equation}
where a $W_n$ action on $E$ is by deck transformation.
\end{cor}
\vskip 2mm
\hskip 120mm $\Box$
\vskip 2mm
\noindent
\begin{rem}
The $W_n$ action on $T^n$ is higly noneffective. The kernel of it is the commutator subgroup of $W_n,$ which
by the Kurosh subgroup theorem, is a free group of rank $1 + (n-2)2^{n-1}.$ 
\end{rem} 
See \cite[ Exercise 3, p. 212]{Cohen} and the proof of Proposition \ref{rem1}.

\section{Cohomologies of $G_n$}
In this part we shall calculate a cohomology of the group $G_n$ with $\Q$ coefficients (Theorem \ref{easy})
and $\F_2$ coefficients (Theorem \ref{noteasy}). 
We shall apply the Leray - Serre spectral sequence of the fibration (\ref{fib}) and
equivalently Lyndon - Hochschild - Serre spectral sequence for the short exact sequence of groups (\ref{maine})
$$
1\to\Z^n\to G_n\stackrel{\nu}\to W_n\to 1.
$$

\subsection{Hilbert-Poincar\'e series}

\begin{Df}
{\em (\cite[p.230]{AT})} Let $M$ be a topological space. 
For a fixed coefficient field $k$, define the Poincar\'e series of $M$ the formal power series	
$$P(x,k) = \Sigma_{i}a_{i}x^{i}$$
where $a_i$ is the dimension of $H^{i}(M,k)$ as a vector space over $k,$ assuming this dimension
is finite for all $i.$
\end{Df}

\begin{theo}\label{easy}
The rational Hilbert-Poincar\'e series of the space 
$$K(G_n,1) = T^n\times_{W_{n}}E$$ is equal to 
\begin{equation}\label{ratio}
P_{n}(x,\Q) = ((1 + x)(1 + \frac{(1 - (-1)^n))}{2}x^{n} + x(\frac{n-2}{2}(1+x)^{n-1} - \frac{n}{2}(1-x)^{n-1})).
\end{equation}
\vskip 1mm
\noindent
In particular, $P_0(x,\Q) = 1, P_1(x,\Q) = x + 1, P_2(x,\Q) =  x^{3} + 1).$
\end{theo}
\vskip 2mm
\noindent
{\bf Proof:}
We start with Lemma.
\noindent
\begin{lm}
For $p > 1, H^p(W_n,\Q) = 0.$ 
\end{lm}
\vskip 1mm
\noindent
{\bf Proof:} 
We have a short exact sequence of groups related to the abelianization
\begin{equation}\label{rational}
1\to \F_{k}\to W_n\to (\Z_2)^n\to 1,
\end{equation}
where $\F_{k}$ is a non abelian free group of a rank $k = 1 + (n-2)2^{n-1}.$
Hence for $q>1, H^{q}(\F_{k},M) = 0$ for any $\F_{k}$-module $M.$ Similar for any $p\geq 1, H^{p}((\Z_2)^{n},N) = 0$ 
for any $(\Z_2)^n$-rational vector space $N.$   
Applying a Leray-Serre spectral sequence to (\ref{rational})
we have for $i\geq 2, H^{i}(W_n,S) = 0.$ Where $S$ is a $W_n-$ rational vector space.
\vskip 5mm \hskip 120mm $\Box$
\vskip 1mm
\noindent
\begin{cor}
For $p > 1, q\geq 0, E_{2}^{p,q} = H^{p}(W_n, H^{q}(\Z^n,\Q)) = 0$ 
and the differentials $d_i = 0$ for $i\geq 2.$ 
Moreover,  $E_{2}^{0,q} \text{and}\hskip 2mm E_{2}^{1,q}, q\geq 0$ are two non trivial columns of the spectral sequence.
\end{cor}
\vskip 2mm
\noindent
The Hilbert-Poincar\'e polynomial (\ref{ratio}) is the sum $f_0 + f_1,$ where
\begin{equation}\label{Pp1}
f_p = x^{p}\Sigma_{i}\text{dim}(E_{2}^{p,i})x^i.
\end{equation}
\vskip 5mm
Let us start to calculate dimensions of $E_{2}^{p,q} = H^{p}(W_n,H^{q}(\Z^n,\Q))$ for $p = 0, 1$ and $q\geq 0.$
We shall use a $W_n$ action on
$H^{q}(\Z^n,\Q) = \Lambda^{q}(\Q^n),$ which follows from an action $W_n$ on $T^n,$ see (\ref{action}).
\vskip 2mm
We introduce sequences $\epsilon\in\{-1,1\}^{n}.$ Denote by $(-1)\epsilon = -\epsilon = (-\epsilon_1,-\epsilon_2,...,-\epsilon_n).$
Moreover, for $A\subset\{1,2,...,n\}$ the sequence $e^{A}$ has $-1$ exactly on the positions from $A.$
Finally let $1 = (1,1,...,1) := e^{\varnothing}$ and $|\epsilon| := \Sigma_{i}\frac{1-\epsilon_i}{2}$ (the number of $-1$ in the sequence).
\vskip 2mm
By $\Q_{\epsilon}$ we shall understand the rational numbers $\Q$ with the structure of a $W_n$-module such that the $k$-th generator of $W_n$ acts
as multipliaction by $\epsilon_k, 1\leq k\leq n.$ 
In this language $H^{1}(\Z^n,\Q)\simeq \Sigma_{i}\Q_{-e^{\{i\}}}$ as $W_{n}$-module. Moreover,
$H^{\ast}(\Z^n,\Q)\simeq \Lambda^{\ast}(\Q^n)$ is a sum of some $\Q_{\epsilon}.$
Let 
$$h^i(\epsilon) = \text{dim}H^{i}(W_n,\Q_{\epsilon}).$$
From the definition

$$h^{0}(\epsilon) =
  \begin{cases}
    1       & \quad \text{if } \epsilon = 1\\
    0  & \quad \text{if } \epsilon\neq 1
  \end{cases}
\text{and}\hskip 2mm
h^{1}(\epsilon) =
  \begin{cases}
    0       & \quad \text{if } \epsilon = 1\\
    |\epsilon|-1  & \quad \text{if } \epsilon\neq 1
  \end{cases}.
$$

Using (\ref{Pp1}) and the above formulas gives us
$$f_i = x^i\sum_{A\subset\{1,2,...,n\}} x^{|A|}h^{i}((-1)^{|A|}e^{A}), i = 1,2.$$

Hence for $n\geq 0$:

$$f_0 = 1 + \frac{1-(-1)^n}{2}x^n$$

$$
f_1 = x^{n+1}(n-1)\frac{1+(-1)^n}{2} +
$$
$$
+ x\hskip 1mm \Sigma_{0<k<n}x^k\binom{n}{k}((k-1)\frac{1+(-1)^k}{2} + (n-k-1)\frac{1-(-1)^k}{2}).
$$
\vskip 5mm
\noindent
\begin{lm}
Formula (\ref{ratio}) from Theorem \ref{easy} is equal to $f_0 + f_1.$
\end{lm}
{\bf Proof}
We shall use two formulas:
$$\Sigma_{0<k<n}\binom{n}{k}x^k = (1+x)^n - 1 - x^n = g(x)$$ and
$$\Sigma_{0<k<n}k\binom{n}{k}x^k = nx(1+x)^{n-1} - nx^n = f(x).$$
On the begining we shall prove that
$$S = \Sigma_{0<k<n}x^k\binom{n}{k}((k-1)\frac{1+(-1)^k}{2} + (n-k-1)\frac{1-(-1)^k}{2})$$ =\footnotemark
\footnotetext{$k(-1)^k +\frac{n-2}{2}-\frac{n}{2}(-1)^k = (k-1)\frac{1+(-1)^k}{2} + (n-k-1)\frac{1-(-1)^k}{2}$}
$$\Sigma_{0<k<n}x^k\binom{n}{k}(k(-1)^k + \frac{n-2}{2} - \frac{n}{2}(-1)^k) =$$
$$\Sigma_{0<k<n}x^k(-1)^k\binom{n}{k} +\frac{n-2}{2}\Sigma_{0<k<n}x^k\binom{n}{k} -\frac{n}{2}\Sigma_{0<k<n}(-1)^kx^k\binom{n}{k} =$$
$$f(-x) + \frac{n-2}{2} -\frac{n}{2}g(-x) = $$
$$\frac{n-2}{2}(1+x)^n - \frac{n}{2}(1+x)(1-x)^{n-1} + (1 - n\frac{1+(-1)^n}{2})x^n +1.$$
Since $$f_1 = x^{n+1}(n-1)\frac{1+(-1)^n}{2} + xS$$ then
$$f_1 = x + \frac{1-(-1)^n}{2}x^{n+1} + (1+x)x(\frac{n-2}{2}(1+x)^{n-1} - \frac{n}{2}(1-x)^{n-1})$$
and $$f_0 + f_1 =$$ 
$$(1+x)(1 + \frac{1-(-1)^n}{2}x^n + x(\frac{n-2}{2}(1+x)^{n-1} - \frac{n}{2}(1-x)^{n-1})).$$
\vskip 5mm
\hskip 120mm $\Box$
\vskip 5mm
\noindent
As a complement to the above results we present an observation about the algebra structure of $H^{\ast}(G_n,\Q).$
\begin{cor} 
	If $x,y\in H^{\ast}(G_n,\Q)$ are such that deg$(x) > 0$ and deg$(y) > 0$ then $xy = 0.$
\end{cor}
\vskip 1mm
\noindent
{\bf Proof:} For $n = 1$ it is obvious. Let us assume $n\geq 2.$
From the proof of Theorem \ref{easy} $H^{1}(G_n,\Q) = 0$ and $H^{s}(G_n,\Q) = 0$ for $s > n+1.$
So, if deg$(x)\geq n$ or deg$(y)\geq n$ then $xy = 0.$ Hence we can assume that deg$(x) < n$ and deg$(y) < n.$
From the proof of Theorem \ref{easy} we know that the appriopriate Serre spectral sequence converging to $H^{\ast}(G_n,\Q)$
has the property:
$$(\star)\hskip 3mm E_{\infty}^{p,q}\neq 0 \Longrightarrow (p,q)\in\{(0,0),(0,n)\}\cup \{(1,i) : i\leq n\}.$$
Let $(F_i)$ denote the filtration of $H^{\ast}(G_n,\Q)$ associated with the spectral sequence. From $(\star)$ and 
deg$(x) < n$ and deg$(y) < n$ it follows that $x,y\in F_1.$ Consequently $xy\in F_2.$ But again from $(\star)$
$F_2 = 0.$
\vskip 5mm
\hskip 120mm $\Box$
\vskip 5mm
A calculation of the cohomology with $\F_2$ coefficients needs different tools.
We shall also apply the Lyndon - Hochschild - Serre spectral sequence for the short exact sequence of groups (\ref{maine})
$$
1\to\Z^n\to G_n\stackrel{\nu}\to W_n\to 1.
$$
In this case we have $E_{3}^{*,*}\neq 0$ and we shall use multiplicative structures. 
\vskip 2mm
\noindent
\begin{theo}\label{noteasy}
\begin{equation}\label{f2}
P_{n}(x,\F_2) = (1 + x)(1 + (n - 1)x(1 + x)^{n-1}).
\end{equation}
\vskip 1mm
\noindent
In particular, $P_0(x,\F_2) = 1, P_1(x,\F_2) = x + 1, P_2(x,\F_2) = 1 + 2x + 2x^2 + x^{3}).$
\end{theo}
\vskip 2mm
\noindent
{\bf Proof:} 
There are canonical isomorphisms over $\F_2$ 
$$E_{2}^{p,q} = H^{p}(W_{n},H^{q}(\Z^n,\F_2))\stackrel{\sim}{\leftarrow}H^{p}(W_n,\F_2)\otimes H^q(\Z^n,\F_2),$$
$$H^{\ast}(\Z_{2}^{n},\F_2) = \F_{2}[z_1,z_2,...,z_n],$$ 
where $z_i\in H^1(\Z_{2}^{n},\F_2) =\hskip 1mm \text{Hom}(\Z_{2}^{n},\F_2)$ is the projection
on the $i$-coordinate, $i=1,2,...,n$ and 
\begin{equation}\label{wn}
H^{\ast}(W_n,\F_2) = H^{\ast}(\Z_{2}^{n},\F_2)/\{z_{i}z_{j}|i\neq j\}.
\end{equation} 
Let $g_1,...,g_n\in H^1(\Z^n,\F_2)$ be a dual basis to $x_{1}^{2},...,x_{n}^{2}.$ 
We shall denote by $\Lambda(g_1,...,g_n)$ the exterior algebra over $\F_2$ generated by 
$g_1,...,g_n$ ($H^{*}(\Z^n,\F_2)$)\footnotemark.
\footnotetext{$\Lambda^{*}(g_1,...,g_n)\simeq H^{*}(\Z^n,\F_2)$}
To begin we shall prove.
\begin{prop}\label{prop2} 
Let $(E_r,d_r)$ be the above spectral sequence and on $E_3$ we use the total grading.
$E_{3}^{0,0} = \langle 1 \rangle, E_{3}^{0,q} = 0, q > 0$ and $E_{3}^{p,q} = 0$ for $p > 2.$ 
Moreover $d_2\neq 0$ and $d_i = 0$ for $i\geq 3.$
\end{prop}
\vskip 1mm
\noindent
{\bf Proof:} 
To prove that $E_{3}^{0,q} = 0, q > 0$ it is enough to show that the kernel of $d_{2}^{0,q}$ is trivial.
By contradiction assume that $\exists\hskip 2mm 0\neq\omega\in E_{2}^{0,q}, q > 0$ and $d_{2}^{0,q}(\omega) = 0.$
Then $\exists\hskip 2mm i,$ s.t. $\omega = g_i\alpha + \beta$ and $\alpha,\beta$ do not depend on $g_i$.
If $d_{2}^{0,q}(\omega) = 0$ then from the properties of the transgression $0 = d_{2}^{0,q}(\omega) = z_{i}\alpha + \gamma,$
where $\gamma$ is a linear combination of elements from the set $\{z_{s}^{2}: s\neq i\}$ with coefficients
in $\Lambda(g_1,...,g_n).$ Hence $\alpha = 0$ a contradiction.
\vskip 1mm
\noindent
From Lemmas \ref{important} and \ref{incli} cycles of the differential $d_2$ are linear combinations of elements
$z_{i}^{k}\omega$ (for some $i,k$) and $\omega\in\Lambda(g_1,g_2,...,g_n)$ does not include $g_i.$
Hence for $k\geq 3\hskip 2mm z_{i}^{k}\omega = d_{2}^{k,s}(z_{i}^{k-2}g_{i}\omega)$ and for
$p\geq 3, E_{3}^{p,q} = 0.$ 
\vskip 1mm
\hskip 120mm $\Box$
\begin{lm}\label{lem5}
In the spectral sequence of the extension (\ref{dihed})
$$0\to\Z\to W_2\to\Z_2\oplus\Z_2\to 0,$$
$\bar{d_{2}}(g) = z_{1}z_{2},$ where $g$ is a generator of $H^{1}(\Z,\F_2).$
\end{lm}
\noindent
{\bf Proof:} From the above  $z_1z_2 \in H^{\ast}(W_2,\F_2)$ is equal to zero. Applying the five-term exact sequence (see \cite[p. 16,57]{AH})
we get the exact sequence
$$H^1(\Z,\F_2)\stackrel{\bar{d_2}}\to H^2(\Z_2\oplus\Z_2,\F_2)\stackrel{f^{\ast}}\to H^2(W_2,\F_2)$$
and Im$\bar{d_2}$ = Ker$f^{\ast}.$ Hence $\bar{d_2}(g) = z_1z_2.$
\vskip 2mm \hskip 120mm $\Box$
\begin{lm}\label{important}
Let $g_1,...,g_n\in H^{1}(\Z^n,\F_2)$ be a dual basis to $x_{1}^{2},...,x_{n}^{2}.$ For the spectral sequence
of the exact sequence of groups
$$0\to\Z^n\to G_n\to W_n\to 0,$$
using naturality and properties of homomorphisms 
(\ref{fii}) $\phi_{i}:G_n\to W_2$ we obtain 
$d_2(g_i) = z_{i}^2.$
\end{lm}
\noindent
{\bf Proof:} We have a commutative diagram
$$\begin{diagram}
\node{1}\arrow{e}\node{\Z^n}\arrow{s,r}{\alpha_{i}}\arrow{e,t}{i}\node{G_{n}}\arrow{s,r}{\phi_{i}}\arrow{e,t}{q_n}\node{W_{n}}\arrow{s,r}{\gamma_{i}}\arrow{e}\node{1}\\
\node{1}\arrow{e}\node{\Z}\arrow{e,t}{i_{1}}\node{W_{2}}\arrow{e,t}{\alpha}\node{\Z_{2}\oplus\Z_{2}}\arrow{e}\node{1}
\end{diagram}.$$
\vskip 1mm
\centerline{Diagram 2}
\vskip 1mm
Here, $\alpha_i$ and $\gamma_i$ are defined by $\phi_i.$ From naturality $d_2\circ\alpha_{i}^{\ast} = \gamma_{i}^{\ast}\circ\bar{d_2},$
(see Diagram 3)
where $\alpha_{i}^{\ast}, \gamma_{i}^{\ast}$ are the induced maps on cohomology. 
\newpage
$$\begin{diagram}\node{\bar{E}_{2}^{0,1}\simeq H^{1}(\Z,\F_2)}\arrow{s,r}{\alpha_{i}^{\ast}}\arrow{e,t}{\bar{d}_2}\node{\operatorname{\bar{E}_{2}^{2,0}\simeq 
H^{2}(\Z_2\oplus\Z_2,\F_2)}}\arrow{s,r}{\gamma_{i}^{\ast}}\\
\node{E_{2}^{0,1}\simeq H^{1}(\Z^{n},\F_{2})}\arrow{e,t}{d_{2}}\node{\operatorname{E_{2}^{2,0}\simeq H^{2}(W_{n},\F_{2})}}
\end{diagram}$$
\vskip 1mm
\centerline{Diagram 3}
\vskip 1mm
\noindent
We have
$$g_i = \alpha_{i}^{\ast}(g).$$
Hence $$d_2(g_i) = d_2(\alpha_{i}^{\ast}(g)) = \gamma_{i}^{\ast}(\bar{d_2}(g)) \stackrel{\text{Lemma}\hskip 1mm\ref{lem5}}{=} \gamma_{i}^{\ast}(z_{1}z_{2}) = \gamma_{i}^{\ast}(z_1)\gamma_{i}^{\ast}(z_2).$$
Moreover, let $\overline{\gamma_{i}}:W_n/[W_n,W_n]\to \Z_{2}\oplus\Z_{2}$ be induced by $\gamma_{i},$ then
$\overline{\gamma_{i}}(\lambda_1, \lambda_2, \cdots ,\lambda_n) = (\Sigma_{s}\lambda_s,\lambda_i).$
Hence, $$\gamma_{i}^{\ast}(z_1) = \Sigma_{s}z_{s}\hskip 2mm \text{and}\hskip 2mm \gamma_{i}^{\ast}(z_2) = z_i.$$
Finally, since in $H^{\ast}(W_n)\hskip 2mm z_{i}z_{j} = 0$ for $i\neq j,$ we get
$$d_2(g_i) = (\Sigma_{s}z_{s})z_{i} = z_{i}^{2}.$$ 
\vskip 2mm \hskip 120mm $\Box$
\vskip 5mm
\noindent

The next observation is the following.
\begin{lm}\label{incli}
Let $k > 0,$ then,
$d_{2}(z_{i}^{k}\omega) = 0$ if and only if $\omega$ does not depend on $g_i.$
\end{lm}
\vskip 1mm
\noindent
{\bf Proof:} Let us write $\omega = g_i\alpha +\beta$ where 
$\alpha,\beta \in \Lambda$ do not depend on $g_i.$
In fact, by definition $d_2(g_i\alpha + \beta) = z_{i}^{2}\alpha + \gamma,$
where $\gamma$ is a linear combination of elements from the set $\{z_{s}^{2} : s\neq i\}$
with coefficients from $\Lambda.$ Hence
, because for $i\neq s, z_{i}z_{s} = 0$ (\ref{wn})
\begin{equation}\label{formulka1}
d_{2}(z_{i}^{k}(g_{i}\alpha + \beta)) = z_{i}^{k}d_{2}(g_{i}\alpha + \beta) = z_{i}^{k}(z_{i}^{2}\alpha + \gamma) = z_{i}^{k+2}\alpha.
\end{equation}
Assume $d_{2}(z_{i}^{k}\omega) = 0.$ We can write $\omega = g_i\alpha + \beta,$ where $\alpha$ and $\beta$ are independent from $g_i.$
From (\ref{formulka1}) $0 = d_2(z_{i}^{k}\omega) = z_{i}^{k}\alpha.$ Hence $\alpha = 0$ and $\omega = \beta$ and so
$\omega$ does not include $g_i.$ Finally, if $\omega$ does not include $g_i$ substituting $\alpha = 0$ and $\omega = \beta$
to the formula (\ref{formulka1}) we get $d_2(z_{i}^{k}\omega) = 0.$
\vskip 2mm
\hskip 120mm $\Box$
\begin{cor}\label{cor1}
Let $k>0$ and $\upsilon = \Sigma_{s}z_{s}^{k}\omega_{s}\in E_{2}^{k,s},$ where $\forall\hskip 1mm s\hskip 1mm \omega_s\in\Lambda_s$ then
$$d_2(\upsilon) = 0 \Longleftrightarrow \forall s\hskip 2mm \omega_s\hskip 2mm \text{does not include}\hskip 2mm g_s.$$ 
\end{cor}
\vskip 1mm
\noindent
{\bf Proof:} ($\Leftarrow$) Follows from the above Lemma \ref{incli}. ($\Rightarrow$) 
For any $i$ if $d_2(\upsilon) = 0$ then also $d_2(z_{i}\upsilon) = 0.$ But $z_i\upsilon = z_{i}^{k+1}\omega_{i}.$
Again, from Lemma \ref{incli} it follows that $\omega_i$ does not include $g_i.$
\vskip 1mm
\hskip 120mm $\Box$

\begin{rem}
Let $Z_{2}^{i,j} =$ker$d_{2}^{i,j}$ and $B_{2}^{i,j} =$Im$d_{2}^{i,j}.$ Let $M$ be an $\F_2$-vector space
and a trivial $G$-module, then 
	$$\dim H^{i}(G,M) = \dim H^{i}(G,\F_2)\dim M.$$ 
	Moreover $\dim H^{i}(\Z^n,\F_2) = \binom{n}{i},$
        $\dim H^{i}(W_n,\F_2) = n$ (for $i > 0$), $\dim E_{2}^{p,q} = n\binom{n}{q}$ (for $p > 0$).
\end{rem}


Summing up the generating function for $H^{\ast}(G_n,\F_2)$ is a sum of three components: $f_0 + f_1 + f_2$ where
\begin{equation}\label{Pp2}
f_p = \Sigma_{i}\text{dim}(E_{3}^{p,i})x^{p+i}.
\end{equation}
\begin{lm}\label{illu}
From the properties of the differentials $d_2$ we have: 
\begin{itemize}
\item[I.] $f_0 = 1;$
\item[II.] $f_1 = nx(1+x)^{n-1};$ 
\item[III.]$f_2 = nx^{2}(1+x)^{n-1} - x((1+x)^{n}-1).$ 
\end{itemize} 
\end{lm}
\vskip 1mm
\noindent
{\bf Proof:}
By an application of the proof of Proposition \ref{prop2} $f_{0} = 1.$

From the above $d_2(z_{i}^{k}\omega) = 0$ if and only if $\omega$ does not include $g_i.$
Hence dim$Z_{2}^{k,s} = n\binom{n-1}{s}$ and $f_1 = nx(1+x)^{n-1},$ c.f. Corollary \ref{cor1}.

For $p = 2$ we have $E_{3}^{2,i} = \text{Ker}d_{2}^{2,i}/\text{Im}d_{2}^{0,i+1}.$
Moreover, for $i > 0, d_{2}^{0,i}$ is a monomorphism. This follows from the proof of Proposition \ref{prop2}.
Hence, dim Im$d_{2}^{0,i}$ = dim$E_{2}^{0,i}$ = $\binom{n}{i}.$
Summing up $f_2 = nx^{2}(1+x)^{n-1} - x((1+x)^{n}-1)$ and the Lemma is proved.
\vskip 2mm
\hskip 120mm $\Box$
\noindent
\vskip 2mm
\begin{exe}
We have dim$Z_{2}^{k,s} = n\binom{n-1}{s}.$ In fact, a basis of $Z_{2}^{k,s}$ is the set $\{z_{i}^{k}\omega\},$
where $1\leq i\leq n$ and $\omega$ is a Grassmann monomial of degree $s$ on elements $\{g_1,...,g_n\}\setminus\{g_i\}.$
\vskip 1mm
\noindent
For $n = 4$ the basis of $Z_{2}^{2,2}$ has 12 elements:
\vskip 1mm
\noindent
$z_{1}^{2}g_{2}g_{3}, z_{1}^{2}g_{2}g_{4},z_{1}^{2}g_{3}g_{4},z_{2}^{2}g_{1}g_{3},z_{2}^{2}g_{1}g_{4},z_{2}^{2}g_{3}g_{4},\\
z_{3}^{2}g_{1}g_{2}, z_{3}^{1}g_{2}g_{4}, z_{3}^{2}g_{2}g_{4}, z_{4}^{2}g_{1}g_{2},z_{4}^{2}g_{1}g_{3},z_{4}^{2}g_{2}g_{3}$
\vskip 1mm
\noindent
the basis of $E_{2}^{0,3}$ has the following 4 elements:
\vskip 1mm
\noindent
$$
y_1 = g_{1}g_{2}g_{3}, y_2 = g_{1}g_{2}g_{4}, y_3 = g_{1}g_{3}g_{4}, y_4 = g_{2}g_{3}g_{4}.
$$
\vskip 1mm
\noindent
Since $d_{2}^{0,3}$ is a monomorphism the basis of $B_{2}^{2,2}$ has four elements:
\vskip 1mm
\noindent
$d_{2}^{0,3}(y_1) = z_{1}^{2}g_{2}g_{3} + z_{2}^{2}g_{1}g_{3} + z_{3}^{2}g_{1}g_{2},$
\vskip 1mm
\noindent
$d_{2}^{0,3}(y_2) = z_{1}^{2}g_{2}g_{4} + z_{2}^{2}g_{1}g_{4} + z_{4}^{2}g_{1}g_{2},$
\vskip 1mm
\noindent
$d_{2}^{0,3}(y_3) = z_{1}^{2}g_{3}g_{4} + z_{3}^{2}g_{1}g_{4} + z_{4}^{2}g_{1}g_{3},$
\vskip 1mm
\noindent
$d_{2}^{0,3}(y_4) = z_{2}^{2}g_{3}g_{4} + z_{3}^{2}g_{2}g_{4} + z_{4}^{2}g_{2}g_{3}.$
\vskip 2mm
\noindent
Hence dim$E_{3}^{2,2} = 12 - 4 = 8.$
\end{exe}
\vskip 5mm
\noindent
Finally
\vskip 5mm
\noindent
$$f_0 + f_1 + f_2 = 1 + nxb + nx^{2}b - xb(1 + x) + x =$$  
$$= x + 1 + nxb(1 + x) - xb(1 + x) = (x + 1)(1 + nxb - xb)=$$
$$= (1 + x)(1 + (n - 1)x(1 + x)^{n-1}).$$
\vskip 2mm
\noindent
Here $b = (1 + x)^{n-1}.$
\vskip 1mm
\noindent
\vskip 5mm
\noindent
Finally, we would like to present some grading of $H^{\ast}(G_n,\F_2).$
We start with a definition.
\begin{Df}\label{grada}
	Define the bigraded algebra $\mathcal{E}^{(n)}$ over $\F_2$ by (a direct sum of vector space):
	$$\mathcal{E}^{(n)} = \mathcal{E}_0 + \mathcal{E}_1 + \mathcal{E}_2$$
where $\mathcal{E}_i$ are given by:
\begin{itemize}
\item $\mathcal{E}_0 = \langle 1\rangle$
\item $\mathcal{E}_{1}$ is spaned (i.e. is a free vector space) by symbols: $z_{i}g_{A}$ where $1\leq i\leq n$ and $A\subset\{1,2,...,n\}, i\notin A$
\item $\mathcal{E}_{2} = \mathcal{E}_{2}'/\mathcal{R}$ where $\mathcal{E}_{2}'$ is spanned by symbols $z_{i}^{2}g_{A}$ with
restrictions as above and $\mathcal{R} =$span$\{r_{A} : A\subset\{1,2,...,n\}, A\neq\emptyset\}$ where 
$r_{A} = \Sigma_{i\in A} z_{i}^{2}g_{A\setminus\{i\}}$
\end{itemize}
Bidegrees are given by:
\vskip 3mm
bideg(1) = (0,0), bideg($z_{i}g_{A}$) = (1,$|A|$), bideg($z_{i}^{2}g_{A}) =$ (2,$|A|$)
\vskip 2mm	
Multiplication is given by:
\begin{itemize}
\item $1$ acts in obvious way;
\item $(z_{i}g_{A})(z_{i}g_{B}) = z_{i}g_{A\cup B}$ if $A\cap B = \emptyset,$
\item All other products are zero.
\end{itemize}
\end{Df}
\vskip 2mm
\noindent
The above definition summarizes explicitly the description of the bigraded algebra structure of $E_3,$ namely
\begin{prop}\label{propcor}
If $(E_r,d_r)$ is the spectral sequence of the short exact sequence (\ref{maine}) then
$$E_3\simeq \mathcal{E}^{(n)}$$
as bigraded algebras.
\end{prop}
\vskip 3mm
\begin{exe}
(Bigraded algebra $H^{\ast}(G_2,\F_2)$)
There are elements: $a,b,A,B,w \in H^{\ast}(G_2,\F_2)$ such that:
\begin{itemize}
	\item $a,b\in H^{1}(G_2,\F_2), A,B \in H^2(G_2,\F_2), w \in H^3(G_2,\F_2);$
\item $aA = bB = w$ and all other products of elements from $\{a,b,A,B,w\}$ are zero.
\item $\{a,b,A,B,w\}$ is a basis of $H^{\ast}(G_2,\F_2).$ 
\end{itemize}
Using Definition \ref{grada} we have:
\vskip 1mm
\noindent
$\mathcal{E}_{0} = \langle 1\rangle,$
\vskip 1mm
\noindent
$\mathcal{E}_{1} = \langle z_{1}g_{\emptyset}, z_{2}g_{\emptyset}, z_{1}g_{\{2\}}, z_{2}g_{\{1\}}\rangle,$
\vskip 1mm
\noindent
$\mathcal{E}_{2}' = \langle z_{1}^{2}g_{\emptyset}, z_{2}^{2}g_{\emptyset}, z_{1}^{2}g_{\{2\}},z_{2}^{2}g_{\{1\}}\rangle,$
\vskip 1mm
\noindent
$\mathcal{R} = \langle z_{1}^{2}g_{\emptyset}, z_{2}^{2}g_{\emptyset}, z_{1}^{2}g_{\{2\}} + z_{2}^{2}g_{\{1\}}\rangle.$
So
$$(1, z_{1}g_{\emptyset},z_{2}g_{\emptyset}, z_{1}g_{\{2\}}, z_{2}g_{\{1\}}, [z_{1}^{2}g_{\{2\}}])$$
is a basis of $\mathcal{E}^{(2)},$ where $[\xi]$ denotes the class of $\xi.$
\vskip 2mm
\noindent
If $(1,a,b,A,B,w)$ are elements of $H^{\ast}(G_2,\F_2)$ which correspond (in this order) to elements of the above basis
then they satisfy the conditions stated above.
\end{exe}
\begin{exe}
Let $X = P_2\bigvee P_2\bigvee S^3$ where $S^3$ is the 3-dimensional sphere and $P_2$ is
the 2-dimensional real projective space. Let $Y = BG_2.$ Then the Poincar\'e polynomials
of $X$ and $Y$ over $\F_2$ and over $\Q$ are the same but $H^{\ast}(X,\F_2)$ and $H^{\ast}(Y,\F_2)$
are not isomorphic as algebras.
\end{exe}
\section{Additional observation}
\begin{prop}
\vskip 1mm
\noindent
1. For $n > 1, (G_{n})_{ab}\simeq \Z_{4}^{n};$
\vskip 1mm
\noindent
2. For $n > 1$ the center of $G_n$ is trivial;
\vskip 1mm
\noindent
3. The Euler characteristic and the first Betti number of  $G_n$ are equal to zero.
\end{prop}
\vskip 1mm
\noindent
{\bf Proof:} 
1. Follows from a direct calculation.
\vskip 1mm
\noindent
2.  From (\ref{maine}) and the fact that the center of $W_n$ is trivial we have an inclusion $Z(G_n)\subset\Z^n.$
Let $v\in Z(G_n)$. From the above $v = \Pi_{i}(x_{i}^{2})^{\alpha_i}$ for some $\alpha_{i}\in\Z.$
Using relations in $G_n,$ we have
$$x_{1}vx_{1}^{-1} = (x_{1}^{2})^{\alpha_{1}}\Pi_{i\geq 2}(x_{i}^{2})^{-\alpha_{i}}.$$
Since $v = x_{1}vx_{1}^{-1},$ then $\alpha_i = 0$ for $i\geq 2.$
Similar $x_{2}vx_{2}^{-1} = v,$ which gives us $\alpha_{1} = 0.$
\vskip 1mm
\noindent
3. From the properties of the Euler characteristic of a fibration (\ref{fib}) 
(see \cite[page 481]{Spanier})
$\chi(K(G_n,1))= \chi(T^n)\chi(E/W_n) = 0\cdot\chi(E/W_n) = 0.$
The conclusion about the Betti number follows directly from Theorem 1.
\vskip 1mm
\vskip 2mm \hskip 120mm $\Box$

\begin{prop}\label{rem1} Let $n\geq 2.$
The short exact sequence of groups (\ref{maine}) defines 
a reprepresentation 
$$h:W_n\to GL(n,\Z), \forall x\in W_n\hskip 2mm h(x)(e_i) = \bar{x}e_i\bar{x}^{-1},$$
where $e_i\in\Z^n$ is the standard basis $i = 1,2,...,n$ and $\nu(\bar{x}) = x.$
However, $K =$ Ker$h \neq 0,$ because $[W_n,W_n]\subset$ Ker$h.$ In particular, 
$K$
is a finitely generated free group of rank $1+ (n-2)2^{2[n/2]-1}.$
\end{prop}
\vskip 1mm
\noindent
{\bf Proof:} From definition of $h$ 
we have an extension  $K\to W_{n}\to\Z_{2}^s,$ where $s = 2[\frac{n}{2}]$ and
$[x$] is the largest integer not exceeding $x.$
Again by the definition of $h$ the commutator subgroup $W'_{n}$ of $W_n$ has index 1 for $n$ even and index 2 for $n$ odd in the group $K.$
Hence $s = 2[\frac{n}{2}].$
Computing (fractional) Euler characteristics we get
(see \cite[Corollary 5.6, p. 245]{B}) from the Euler characteristic formula 
$e(K) = e(W_n)/e(\Z_{2}^{s}) = e(W_n)2^s = (1 - \frac{n}{2})2^s.$ Which gives
the announced rank.
\vskip 1mm
\noindent
Analogously $e(W'_{n}) = (1 - \frac{n}{2})2^n$ and $e(K/W'_{n}) = 2^{s-n}.$
\vskip 1mm
\hskip 120mm $\Box$

\begin{rem}
Let $n$ be an even number then
$$P_{n}(x,\Q) = P_{n}(x,\F_2)\hskip 5mm \text{mod}\hskip 2mm 2.$$
\end{rem}


\end{document}